\documentclass[12pt,a4paper]{article}
\usepackage[dvips]{graphicx,color} 
\usepackage[T1]{fontenc}  
\usepackage{times}    
\usepackage{enumerate,amsmath,amsthm,mathrsfs,citesort}
\usepackage{amsfonts}
\usepackage{amssymb}
\usepackage[a4paper,left=2.5cm,right=2cm,top=2.5cm,bottom=2.5cm]{geometry}
\usepackage{amstext} 
\usepackage{array}   
\usepackage[textsize=scriptsize,textwidth=4cm,linecolor=cyan, backgroundcolor=cyan]{todonotes}

\newtheoremstyle{bthm}{\baselineskip}{\baselineskip}{\slshape}{}{\bfseries}{}{  }{}
\newtheoremstyle{bex}{\baselineskip}{\baselineskip}{}{}{\itshape}{}{  }{}
\theoremstyle{bthm}
\newtheorem{theorem}{Theorem}[section]
\newtheorem{corollary}[theorem]{Corollary}

\newtheorem{proposition}[theorem]{Proposition}

\newtheorem{fact}[theorem]{Fact}

\theoremstyle{bex}

\begin{document}
\begin{titlepage}
\title{$\chi$-binding functions for some classes of $(P_3\cup P_2)$-free graphs}
\author{Athmakoori Prashant$^{1}$, P. Francis$^{2}$ and S. Francis Raj$^{1}$}
\date{{\footnotesize$^{1}$ Department of Mathematics, Pondicherry University, Puducherry-605014, India}\\
{\footnotesize$^{2}$ Department of Computer Science, Indian Institute of Technology, Palakkad 678557, India}\\
{\footnotesize 11994prashant@gmail.com, pfrancis@iitpkd.ac.in, francisraj\_s@pondiuni.ac.in }}
\maketitle
\renewcommand{\baselinestretch}{1.3}\normalsize
\begin{abstract}
The class of $2K_2$-free graphs have been well studied in various contexts in the past.
It is known that the class of $\{2K_2,2K_1+K_p\}$-free graphs and $\{2K_2,(K_1\cup K_2)+K_p\}$-free graphs admits a linear $\chi$-binding function.
In this paper, we study the classes of $(P_3\cup P_2)$-free graphs which is a superclass of $2K_2$-free graphs.
We show that $\{P_3\cup P_2,2K_1+K_p\}$-free graphs and $\{P_3\cup P_2,(K_1\cup K_2)+K_p\}$-free graphs also admits linear $\chi$-binding functions.
In addition, we give tight chromatic bounds for $\{P_3\cup P_2,HVN\}$-free graphs and $\{P_3\cup P_2,diamond\}$-free graphs and it can be seen that the latter  is an improvement of the existing bound given by A. P. Bharathi and S. A. Choudum [Colouring of $(P_3\cup P_2)$-free graphs, Graphs and Combinatorics 34 (2018), 97-107].
\end{abstract}
\noindent
\textbf{Key Words:} Chromatic number, $\chi$-binding function, $(P_3\cup P_2)$-free graphs and Perfect graphs. \\
\textbf{2000 AMS Subject Classification:} 05C15, 05C75

\section{Introduction}\label{intro}
All graphs considered in this paper are simple, finite and undirected.
Let $G$ be a graph with vertex set $V(G)$ and edge set $E(G)$.
For any positive integer $k$, a \emph{proper $k$-coloring} of a graph $G$ is a mapping $c$ : $V(G)\rightarrow\{1,2,\ldots,k\}$ such that adjacent vertices receive distinct colors. 
If a graph $G$ admits a proper $k$-coloring, then $G$ is said to be \emph{$k$-colorable}.
The \emph{chromatic number}, $\chi(G)$, of a graph $G$ is the smallest $k$ such that $G$ is $k$-colorable.
Let $P_n, C_n$ and $ K_n$ respectively denote the path, the cycle and the complete graph on $n$ vertices.
For $S,T\subseteq V(G)$, let $N_T(S) = N(S)\cap T$ (where $N(S)$ denotes the set of all neighbors of $S$ in $G$),
let $\langle S\rangle$ denote the subgraph induced by $S$ in $G$ and let $[S,T]$ denote the set of all edges with one end in $S$ and the other end in $T$.
If every vertex in $S$ is adjacent with every vertex in $T$, then $ [S, T ]$ is said to be complete.
For any graph $G$, let $\overline{G}$ denote the complement of $G$.

Let $\mathcal{F}$ be a family of graphs.
We say that $G$ is \emph{$\mathcal{F}$-free} if it does not contain any induced subgraph which is isomorphic to a graph in $\mathcal{F}$.
For a fixed graph $H$, let us denote the family of $H$-free graphs by $\mathcal{G}(H)$.
If $G_1$ and $G_2$ are two vertex-disjoint graphs, then their \emph{union} $G_1\cup G_2$ is the graph with vertex set $V(G_1\cup G_2)=V(G_1)\cup V(G_2)$ and edge set $E(G_1\cup G_2)=E(G_1)\cup E(G_2)$.
Similarly, the \emph{join} of $G_1$ and $G_2$, denoted by $G_1+G_2$, is the graph whose vertex set $V(G_1+G_2) = V(G_1)\cup V(G_2)$ and the edge set $E(G_1+G_2) = E(G_1)\cup E(G_2)\cup\{xy: x\in V(G_1),\ y\in V(G_2)\}$. Let $\omega(G)$ and $\alpha(G)$ denote the \emph{clique number} and  \emph{independence number} of a graph $G$ respectively.
When there is no ambiguity, $\omega(G)$ will be denoted by $\omega$.

A graph $G$ is said to be \emph{perfect} if $\chi(H)=\omega(H)$, for every induced subgraph $H$ of $G$.
One of the well-known conjectures on graph colorings was given by C. Berge in \cite {berge1961farbung} and later proved by
M. Chudnovsky et al., in \cite{chudnovsky2006strong}.

\begin{theorem}(\cite{chudnovsky2006strong})\label{SPGT}
\textbf{The Strong Perfect Graph Theorem}.
A graph is perfect if and only if it does not contain $C_{2k+1}$ or $\overline{C}_{2k+1}$ as induced subgraphs, for any $k\geq 2$.
\end{theorem}

From the Strong Perfect Graph Theorem we see that, $\chi(G)=\omega(G)$ for $\{C_5,C_7,C_9,\ldots,\linebreak\overline{C_5},\overline{C_7},\overline{C_9},\ldots\}$-free graph $G$.
Although, we cannot anticipate the same for other graph classes described in terms of some other forbidden induced subgraphs, the challenge of determining an upper bound on the chromatic number of graphs in terms of their clique number is both fascinating and difficult.
In order to determine an upper bound for the chromatic number of a graph in terms of their clique number, the concept of $\chi$-binding functions was introduced by A. Gy{\'a}rf{\'a}s in \cite{gyarfas1987problems}.
A hereditary class $\mathcal{G}$ of graphs is said to be \emph{$\chi$-bounded} \cite{gyarfas1987problems} if there is a function $f$ (called a $\chi$-binding function) such that $\chi(G)\leq f(\omega(G))$, for every $G\in \mathcal{G}$.
We say that the $\chi$-binding function $f$ is \emph{special linear} if $f(x)= x+c$, where $c$ is a constant.
There has been extensive research done on $\chi$-binding functions for various graph classes.
See for instance, \cite{skarthick2018chromatic,cameron2021optimal,bharathi2018colouring,brause2019chromatic,schiermeyer2019polynomial,karthick2016vizing,randerath2004vertex,gyarfas1987problems,karthick2018chromatic}.
In \cite{gyarfas1987problems}, A. Gy{\'a}rf{\'a}s posed a lot of questions.

The family of $2K_2$-free graphs have been well studied in various contexts: domination (El-Zahar and P. Erd\H{o}s \cite{el1985existence}), size (F. R. K. Chung et al., \cite{chung1990maximum}), vertex coloring (S. Wagon \cite{wagon1980bound}, S. Gaspers and S. Huang \cite{gaspers20192p_2}, A. Gy{\'a}rf{\'a}s \cite{gyarfas1987problems}),
edge coloring (P. Erd\H{o}s \cite{erdos1985problems}) and algorithmic complexity (Z. Bl{\'a}zsik et al., \cite{blazsik1993graphs}).
One of the earliest results was by S. Wagon \cite{wagon1980bound}, in which he showed that the $\chi$-binding function for $2K_2$-free graphs is $\binom{\omega+1}{2}$. Motivated by this A. Gy{\'a}rf{\'a}s in \cite{gyarfas1987problems} posed a problem which asks for the order of magnitude of the smallest $\chi$-binding function for $\mathcal{G}(2K_2)$.
In  \cite{brause2019chromatic}, C. Brause et al., showed that the class of $(2K_2,H)$-free graphs does not admit a linear $\chi$-binding function when $\alpha(H)\geq 3$.
In \cite{karthick2018chromatic}, T. Karthick and S. Mishra focused on finding subfamilies of $2K_2$-free graphs which admit special linear $\chi$-binding functions and proved that the families of $\{2K_2, H\}$-free graphs, where $H\in \{HVN, diamond, K_1+P_4, K_1 + C_4, \overline{P_5}, \overline{P_2 \cup P_3}, K_5-e\}$ admit special linear $\chi$-binding functions.
In particular, they proved that $\{2K_2, K_5-e\}$-free graphs and $\{2K_2, K_1+C_4\}$-free graphs admit the $\chi$-binding functions $\omega(G)+4$ and $\omega(G)+5$ respectively. These bounds were later improved  by Athmakoori Prashant et al., in \cite{prashant2021chromatic,caldam2021chromatic} to $\{4\  (\text{for~} \omega(G)\leq 3), 6\  (\text{for~}\omega(G)=4)$ and $\omega(G) (\text{for~}\omega(G)\geq 5)\}$ and $\omega(G) +1$ respectively. Also, they showed that the class of $\{2K_2,K_2+P_4\}$-free graphs admit a special linear $\chi$-binding function $\omega(G)+2$.
In  \cite{brause2019chromatic}, C. Brause et al., improved the $\chi$-binding function  for $\{2K_2,K_1+P_4\}$-free graphs to $\max\{3, \omega(G)\}$.
Also they proved that for $s\neq 1$ or $\omega(G)\neq 2$, the class of $\{2K_2, (K_1\cup K_2)+K_s\}$-free graphs with $\omega(G)\geq 2s$ is perfect and for  $r\geq 1$, the class of $\{2K_2, 2K_1+K_r\}$-free graphs with $\omega(G)\geq 2r$  is perfect.
Clearly when $s=2$ and $r=3$, $(K_1\cup K_2)+K_s\cong HVN$ and $2K_1+K_r\cong K_5-e$ which implies that the class of $\{2K_2,HVN\}$-free graphs and $\{2K_2,K_5-e\}$-free graphs are perfect for $\omega(G)\geq 4$ and $\omega(G)\geq 6$ respectively which improved the bounds given in \cite{karthick2018chromatic}.

\begin{figure}
  \centering
  \includegraphics[width=16cm]{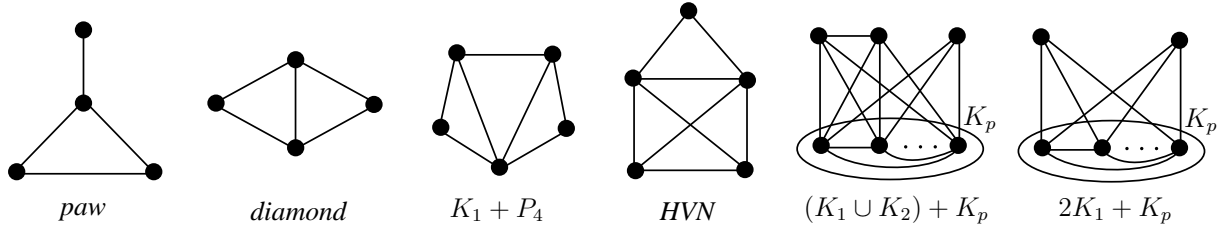}
  \caption{Some special graphs}\label{freegraphs}
\end{figure}

Motivated by C. Brause et al., and their work on $2K_2$-free graphs in \cite{brause2019chromatic}, we started looking at $(P_3\cup P_2)$-free graphs which is a superclass of $2K_2$-free graphs.
In \cite{bharathi2018colouring}, A. P. Bharathi and S. A. Choudum obtained a $O(\omega^3)$ upper bound for the chromatic number of $(P_3\cup P_2)$-free graphs and obtained sharper bounds for $\{P_3\cup P_2,diamond\})$-free graphs.
Also, T. Karthick and S. Mishra in \cite{skarthick2018chromatic} showed that if $G$ is $\{P_3\cup P_2,diamond,K_4\}$-free, then $\chi(G)=6$.
The subclass of $diamond$-free graphs have been studied in the past. See for instance \cite{arbib2002p5,brandstadt2004p5,cameron2021optimal,chiarelli2021strong,choudum2010first,karthick2016vizing,kloks2009even,tucker1987coloring}.
In this paper, by using structural results, we obtain linear $\chi$-binding functions for the class of $\{P_3\cup P_2, (K_1\cup K_2)+K_p\}$-free graphs and $\{P_3\cup P_2, 2K_1+K_p\}$-free graphs.
In addition, for $\omega(G)\geq 3p-1$ we show that the class of $\{P_3\cup P_2, (K_1\cup K_2)+K_p\}$-free graphs admits a special linear $\chi$-binding function $f(x)=\omega(G) +p-1$ and $\{P_3\cup P_2, 2K_1+K_p\}$-free graphs are perfect. We also construct an example of a $\{P_3\cup P_2, (K_1\cup K_2)+K_p\}$-free graph with $\omega(G)\geq 3p-1$ such that $\chi(G)\geq \omega(G)+\lceil\frac{p-1}{2}\rceil$.
In addition, we give a tight $\chi$-binding function for $\{P_3\cup P_2,HVN\}$-free graphs and for $\{P_3\cup P_2,diamond\}$-free graphs.
This bound for $\{P_3\cup P_2, diamond\}$-free graphs turns out to be an improvement of the existing bound obtained by A. P. Bharathi and S. A. Choudum in \cite{bharathi2018colouring}.

Some graphs that are considered as forbidden induced subgraphs in this paper are given in Figure \ref{freegraphs}.
Notations and terminologies not mentioned here are as in \cite{west2005introduction}.
\section{Preliminaries}

Throughout this paper, we use a particular partition of the vertex set of a graph $G$ as defined initially by S. Wagon in \cite{wagon1980bound} and later improved by A. P. Bharathi and S. A. Choudum  in \cite{bharathi2018colouring} as follows.
Let $A=\{v_1,v_2,\ldots,v_{\omega}\}$ be a maximum clique of $G$.
Let us define the \emph{lexicographic ordering} on the set $L=\{(i, j): 1 \leq i < j \leq \omega\}$ in the following way.
For two distinct elements $(i_1,j_1),(i_2,j_2)\in L$, we say that $(i_1,j_1)$ precedes $(i_2,j_2)$, denoted by $(i_1,j_1)<_L(i_2,j_2)$ if either $i_1<i_2$ or $i_1=i_2$ and $j_1<j_2$.
For every $(i,j)\in L$, let $C_{i,j}=\left\{v\in V(G)\backslash A:v\notin N(v_i)\cup N(v_j)\right\}\backslash\left\{\mathop\cup\limits_{(i',j')<_L(i,j)} C_{i',j'}\right\}$.
Note that, for any $k\in \{1,2,\ldots,j\}\backslash\{i,j\}$, $[v_k,C_{i,j}]$ is complete.
Hence $\omega(\langle C_{i,j}\rangle)\leq\omega(G)-j+2$.

For $1\leq k\leq \omega$, let us define $I_k=\{v\in V(G)\backslash A: v\in N(v_i), \text{\ for\ every}\ i\in\{1,2,\ldots,\omega\}\backslash\{k\}\} $.
Since $A$ is a maximum clique, for $1\leq k\leq \omega$, $I_k$ is an independent set and for any $x\in I_k$, $xv_k\notin E(G)$.
Clearly, each vertex in $V(G)\backslash A$ is non-adjacent to at least one vertex in $A$. 
Hence those vertices will be contained either in $I_k$ for some $k\in\{1,2,\ldots,\omega\}$, or in $C_{i,j}$ for some $(i,j)\in L$.
Thus $V(G)=A\cup\left(\mathop\cup\limits_{k=1}^{\omega}I_k\right)\cup\left(\mathop\cup\limits_{(i,j)\in L}C_{i,j}\right)$.
Sometimes, we use the partition  $V(G)=V_1\cup V_2$, where $V_1=\mathop\cup\limits_{1\leq k\leq\omega}(\{v_k\}\cup I_k)=\mathop\cup\limits_{1\leq k\leq\omega}U_k$ and $V_2=\mathop\cup\limits_{(i,j)\in L}C_{i,j}$.

Let us recall a  result on $\left(P_3\cup P_2\right)$-free graphs given by A. P. Bharathi and S. A. Choudum in \cite{bharathi2018colouring}.

\begin{theorem}(\cite{bharathi2018colouring}) \label{P3P2chi}
If a graph $G$ is $(P_3\cup P_2)$-free, then $\chi(G)\leq \frac{\omega(G)\left(\omega(G)+1\right)\left(\omega(G)+2\right)}{6}$.
\end{theorem}

Without much difficulty one can make the following observations on $\left(P_3\cup P_2\right)$-free graphs.
\begin{fact}\label{true} Let $G$ be a  $\left(P_3\cup P_2\right)$-free graph. For $(i,j)\in L$, 
the following holds.
\begin{enumerate}[(i)]
\setlength\itemsep{-1pt}
\item \label{CijClique} Each $\langle C_{i,j}\rangle$ is a disjoint union of cliques, that is, $\langle C_{i,j}\rangle$ is $P_3$-free.
\item \label{joinvertex} For every integer $s\in \{1,2,\ldots,j\}\backslash\{i,j\}$, $N(v_s)\supseteq\{C_{i,j}\cup A\cup(\mathop\cup\limits_{k=1}^{\omega}I_k)\}\backslash\{v_s\cup I_s\}$.
\end{enumerate}
\end{fact}
\begin{proof}
\begin{enumerate}[(i)]\setlength\itemsep{-1pt}
\item  Suppose a $\langle C_{i,j}\rangle$ contains a $P_3$ say $P$, then $\langle \{V(P),v_i,v_j\}\rangle\cong P_3\cup P_2$, contradiction.
\item Follows immediately from the definitions of $C_{i,j}, A$ and $I_k$, $1\leq k\leq \omega(G)$.
\vskip-.6cm
\end{enumerate}
\end{proof}

\begin{figure}
  \centering
  \includegraphics[width=4cm]{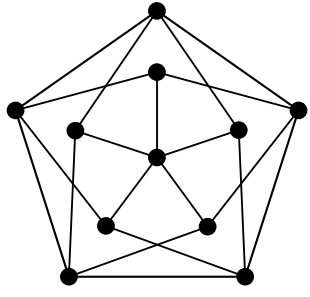}
\caption{Gr\"{o}tzsch graph $\mu(C_5)$}
\label{gro}
\end{figure}

When $\omega(G)=2$, by Theorem \ref{P3P2chi} we see that if $G$ is $(P_3\cup P_2)$-free, then $\chi(G)\leq4$.
If we consider the Gr\"{o}tzsch graph $\mu(C_5)$ in Figure \ref{gro}, we can observe that $\left< V(\mu(C_5))\backslash N(V(P_3))\right>$, for any $P_3$ in $\mu(C_5)$ is either $K_1$, $2K_1$ or $3K_1$ and hence $\mu(C_5)$ is $(P_3\cup P_2)$-free.
Also $\omega(\mu(C_5))=2$ and $\chi(\mu(C_5))=4$.
Thus when $\omega(G)=2$, the bound given in Theorem \ref{P3P2chi} is tight. Therefore while trying to improve the bound given in Theorem \ref{P3P2chi} for $(P_3\cup P_2)$-free graphs we shall consider $\omega\geq 3$.

\section{$\{P_3\cup P_2, (K_1\cup K_2)+K_p\}$-free graphs}

Let us start Section 3 by observing that when $p=0$, $((K_1\cup K_2)+K_p)$-free graphs will be $P_4$-free and hence perfect.
Therefore throughout this section we will assume that $p\geq 1$.
Let us make a few more observations on $((K_1\cup K_2)+K_p)$-free graphs. 

\begin{proposition}\label{k1k2kpprop}
Let $G$ be a $((K_1\cup K_2)+K_p)$-free graph with $\omega(G)\geq p+2$, $p\geq 1$.
Then $G$ satisfies the following.
\begin{enumerate}[(i)]
\setlength\itemsep{-1pt}
\item \label{completeIkIl} For $k,\ell\in\{1,2,\ldots,\omega(G)\}$, $[I_k,I_{\ell}]$ is complete. Thus, $\langle V_1\rangle$ is a complete multipartite graph with  $U_{k}=\{v_{k}\}\cup I_{k}$, $1\leq k\leq \omega(G)$ as its partitions.

\item \label{emptycij} For $j\geq p+2$ and $1\leq i<j$, $C_{i,j}=\emptyset$.

\item \label{p-1neighbors} For $x\in V_2$, $x$ has neighbors in  at most $(p-1)$ $U_\ell$'s where $\ell\in\{1,2,\ldots,\omega(G)\}$.

\end{enumerate}
\end{proposition}
\begin{proof}
Let $G$ be a $((K_1\cup K_2)+K_p)$-free graph with $\omega(G)\geq p+2$, $p\geq 1$.

(\ref{completeIkIl}) For some $k,\ell\in\{1,2,\ldots,\omega\}$, suppose $[I_k,I_{\ell}]$ is not complete, then there exist non-adjacent vertices $a\in I_k$ and $b\in I_{\ell}$ . Since $\omega(G)\geq p+2$, there are integers $q_1, q_2,\ldots, q_p$ in $\{1,2,\ldots, \omega\}\backslash\{k,\ell\}$ such that  $\langle\{a,b,v_k,v_{q_1},\ldots, v_{q_p}\}\rangle\cong (K_1\cup K_2)+K_p$, a contradiction.
Thus, $[I_k,I_{\ell}]$ is complete for any $k,\ell\in\{1,2,\ldots,\omega\}$ and thereby $\langle V_1\rangle$ is a complete multipartite graph with $U_{k}=\{v_{k}\}\cup I_{k}$ as its partitions.

\noindent (\ref{emptycij})  Suppose $C_{i,j}\neq\emptyset$, for some $j\geq p+2$ and $1\leq i<j$, let $a\in C_{i,j}$. If $S\subseteq\{v_1,v_2,\ldots, v_{j-1}\}\backslash\{v_i\}$ such that $|S|=p$, by the definition of $C_{i,j}$, we see that
$\langle\{a,v_i,v_j,S\}\rangle\cong (K_1\cup K_2)+K_p$, a contradiction. Thus, $C_{i,j}=\emptyset$ for all $j\geq p+2$.

\noindent (\ref{p-1neighbors})   Let $x\in V_2$, that is, $x\in C_{i,j}$ for some  $(i,j)\in L$. Suppose there exist $k (\geq p)$ distinct integers $s_1,s_2,\ldots, s_{k}$ in $\{1,2,\ldots, \omega\}$ such that $N(x)\cap U_{s_{\ell}}\neq \emptyset$, $1\leq \ell\leq k$, then we can find $u_{s_{\ell}}\in U_{s_{\ell}}$ where $xu_{s_{\ell}}\in E(G)$.
Let $S=\{u_{s_1},u_{s_2},\ldots, u_{s_k}\}\backslash\{u_i,u_j\}$.
If $k\geq\omega(G)-1$, then $|S|\geq k-2\geq \omega(G)-3 \geq p-1$.  
Suppose $|S|=p-1$, we have $|S|=k-2$ and thus $u_i,u_j\in \{u_{s_1},u_{s_2},\ldots, u_{s_k}\}$ and we get $\langle\{v_i,x,u_i,S\cup\{u_j\}\}\rangle\cong (K_1\cup K_2)+K_p$, a contradiction.
Suppose $|S|\geq p$,  let $S'\subseteq S$ such that $|S'|=p$. Then  $\langle\{x,v_i,v_j,S'\}\rangle\cong (K_1\cup K_2)+K_p$, a contradiction.
If $k\leq\omega(G)-2$, then there exist distinct integers $q,r$ such that  $[x,U_q\cup U_r]=\emptyset$ and thus $\langle\{x,v_q,v_r,u_{s_1},\ldots, u_{s_p}\}\rangle\cong (K_1\cup K_2)+K_p$, a contradiction. Thus (\ref{p-1neighbors}) holds.
\end{proof}
\begin{proposition}\label{p3p2k1k2kpprop}
Let $G$ be a $\{P_3\cup P_2,(K_1\cup K_2)+K_p\}$-free graph with $\omega(G)\geq p+2$, $p\geq 1$.
Then $G$ satisfies the following.
\begin{enumerate}[(i)]
\setlength\itemsep{-1pt}
\item \label{omegaCkj}For $(i,j)\in L$, if $\omega (\langle C_{i,j}\rangle)\geq p-j+4$, then $\omega (\langle C_{k,j}\rangle)\leq1$ for $k\neq i$ and $1\leq k \leq j-1$.

\item \label{p3p2k1k2kpobs}
If $\omega(G)=(p+2+k)$, $k\geq 0$, then $\left\langle \left(\mathop\cup \limits_{j=\max\{2,p+1-\lfloor\frac{k}{2}\rfloor\}}^{p+1} \left(\mathop\cup \limits_{i=1}^{j-1} C_{i,j}\right)\right)\right\rangle$ is $P_3$-free.
\end{enumerate}
\end{proposition}

\begin{proof}
Let $G$ be a $\{P_3\cup P_2,(K_1\cup K_2)+K_p\}$-free graph with $\omega(G)\geq p+2$, $p\geq 1$.

\noindent (\ref{omegaCkj})
For $j\geq p+2$, by (\ref{emptycij}) of Proposition \ref{k1k2kpprop} we see that $C_{i,j}=\emptyset$.
Thus $\omega(\langle C_{i,j}\rangle)=0$ for every $i\in\{1,2,\ldots,j-1\}$.
Next for $j\leq p+1$, let $\omega (\langle C_{i,j}\rangle)\geq p-j+4$.
Suppose  $\omega (\langle C_{k,j}\rangle)\geq2$ for some $k\neq i$ and $1\leq k \leq j-1$, then there exists vertices $a,b\in C_{k,j}$ such that $ab\in E(G)$. Let $S$ be the vertices of a maximum clique in $C_{i,j}$. Then, $|S|=\omega (\langle C_{i,j}\rangle)\geq p-j+4$.
First, we show that $[a,S]$ is complete. On the contrary, let $x\in S$ such that $ax\notin E(G)$.
If there exists a vertex $y\in S\backslash \{x\}$ such that $ay\notin E(G)$ then $\langle\{a,v_i,v_j,x,y\}\rangle\cong P_3\cup P_2$, a contradiction.
Therefore $[a,S\backslash \{x\}]$ is complete and hence for any $S'\subseteq S\backslash{\{x\}}$ with  $|S'|=p-j+3$, $\langle\{a,x,v_k,\{S'\cup\{v_1,v_2,\ldots,v_{j-1}\}\backslash\{v_i,v_k\}\}\}\rangle\cong (K_1\cup K_2)+K_p$, a contradiction.
Hence, $[a,S]$ is complete.
Similarly, we can prove that $[b,S]$ is also complete. Finally, if $S''\subseteq S$ such that $|S''|=p-j+4$, then  we get that $\langle\{v_k,a,b,\{S''\cup\{v_1,v_2,\ldots,v_{j-2}\} \backslash\{v_i,v_k\}\}\}\rangle\cong (K_1\cup K_2)+K_p$, a contradiction.
Hence, $\omega(\langle C_{k,j}\rangle)\leq 1$.

\noindent (\ref{p3p2k1k2kpobs})
Let $\omega(G)=(p+2+k)$, $k\geq 0$.  On the contrary, let $P=\langle\{a,b,c\}\rangle$ be a $P_3$ in \linebreak $ \left\langle \left(\mathop\cup \limits_{j=\max\{2,p+1-\lfloor\frac{k}{2}\rfloor\}}^{p+1} \left(\mathop\cup \limits_{i=1}^{j-1} C_{i,j}\right)\right)\right\rangle$, such that   $a\in C_{i_1,j_1}$, $b\in C_{i_2,j_2}$ and $c\in C_{i_3,j_3}$ where    $\max \{2,p+1-\left\lfloor\frac{k}{2}\right\rfloor\}\leq j_1,j_2,j_3 \leq p+1$, and $ i_{\ell}<j_{\ell} ,1\leq \ell\leq 3$.
Let $j'=\min\{j_1,j_2,j_3\}$. Clearly, $j'\geq \max \{2,p+1-\left\lfloor\frac{k}{2}\right\rfloor\}$.
For $1\leq \ell\leq3$, any vertex $x\in C_{i_{\ell},j_{\ell}}$ has either $(j'-2)$ or $(j'-1)$ neighbors in $\{v_1,v_2,\ldots,v_{j'-1}\}$.
By (\ref{p-1neighbors}) of Proposition \ref{k1k2kpprop}, we have $|N_A (x)|\leq p-1$ and hence the vertex $x$ has at most either $(p+1-j')$ or $(p-j')$ neighbors in $A\backslash\{v_1,v_2,\ldots,v_{j'-1}\}$.
Thus, $|N_A(V(P))|\leq (j'-1)+3(p+1-j')=3p-2j'+2$.
Hence there are at least $\omega (G)-|N_A(V(P))|$ vertices in $A$ which are non neighbors of any of the vertex in $P$. Now,
              $\omega (G)-|N_A(V(P))|\geq (p+2+k)-(3p-2j'+2)$
               $\geq (-2p+k)+2(\max\{2,p+1-\left\lfloor\frac{k}{2}\right\rfloor\})\geq2$.
Hence, there exists vertices $v_r,v_s\in A$ such that $[\{v_r,v_s\},V(P)]=\emptyset$ and thereby we get
$\langle\{V(P),v_r,v_s\}\rangle\cong P_3\cup P_2$, a contradiction.
\end{proof}

As a consequence of Proposition \ref{k1k2kpprop}, we obtain Corollary \ref{k1k2kpcor1} which is a result due to S. Olariu in \cite{olariu1988paw}.
\begin{corollary}\cite{olariu1988paw}\label{k1k2kpcor1}
Let $G$ be a connected graph. Then $G$ is $paw$-free graph if and only if $G$ is either $K_3$-free or complete multipartite.
\end{corollary}
\begin{proof}
For $\omega(G)=1$ or $2$ there is nothing to prove.
Let $\omega(G)\geq 3$.
If $p=1$, then $((K_1\cup K_2)+K_p)\cong paw$. By Proposition \ref{k1k2kpprop}, $V(G)=V_1\cup V_2$ where $V_1$ is a complete multipartite graph, $V_2=C_{1,2}$ and $[C_{1,2},V_1]=\emptyset$.
Since $G$ is connected, we have $C_{1,2}=\emptyset$ and thus $G$  is a complete multipartite graph.
Conversely, if $G$ is either $K_3$-free or a complete multipartite graph, then we can observe that $G$ is $paw$-free.
\end{proof}

Now,  for $p\geq 1$ and $\omega\geq \max\{3,3p-1\}$, let us determine the structural characterization and the chromatic number of  $\{P_3\cup P_2,(K_1\cup  K_2)+K_p\}$-free graphs.

\begin{theorem}\label{k1k2kp}
Let $p$ be a positive integer and $G$ be a $\{P_3\cup P_2, (K_1\cup K_2)+K_p\}$-free graph with $V(G)=V_1\cup V_2$. If $\omega(G)\geq \max\{3,3p-1\}$, then (i) $\langle V_1\rangle$ is a complete multipartite graph with partition $U_1,U_2, \ldots, U_{\omega}$, (ii) $\langle V_2\rangle$ is $P_3$-free graph and
(iii)  $\chi(G)\leq \omega(G)+p-1$.
\end{theorem}
\begin{proof}
Let $p\geq 1$ and $G$ be a $\{P_3\cup P_2, (K_1\cup K_2)+K_p\}$-free graph with $\omega(G)\geq \max\{3,3p-1\}$.
By (\ref{completeIkIl}) of Proposition \ref{k1k2kpprop}, we see that $\langle V_1\rangle$ is a complete multipartite graph with the partition $U_{k}=\{v_{k}\}\cup I_{k}$, $1\leq k\leq \omega$.
By (\ref{CijClique}) of Fact \ref{true}, each $\langle C_{i,j}\rangle$ is $P_3$-free, for every $(i,j)\in L$.
Now, we show that $\left\langle V_2\right\rangle$ is $P_3$-free.
Suppose there exist vertices $a,b,c\in V_2$ such that $\langle\{a,b,c\}\rangle\cong P_3$.
Since $p\geq 1$ and $\omega(G)\geq3$, we have $\omega(G)\geq 3p-1\geq p+2$, and hence by (\ref{p-1neighbors}) of Proposition \ref{k1k2kpprop}, we have $|N_A(\{a,b,c\})|\leq 3(p-1)$. Therefore there exist distinct integers $s,t\in \{1,2,\ldots,\omega\}$ such that $[\{a,b,c\}, \{v_{s}, v_{t}\}]= \emptyset$ and thus $\langle\{a,b,c,v_s,v_t\}\rangle\cong P_3\cup P_2$, a contradiction. Hence, $\langle V_2\rangle$ is $P_3$-free.
Now, let us exhibit an $(\omega+p-1)$-coloring for $G$ using $\{1,2,\ldots, \omega+p-1\}$ colors.
For $1\leq k \leq \omega$, let us give the color $k$ to the vertices of $U_k$.
Let $H$ be a component in $\left\langle V_2\right\rangle$.
Clearly, each vertex in $H$ is adjacent to at most $p-1$ colors given to the vertices of $V_1$ and is adjacent to $\omega(H)-1$ vertices of $H$.
Since $\omega(H)\leq \omega(G)$, each vertex is adjacent to at most $\omega(G)+p-2$ colors when coloring the vertices of $H$. Hence, there is a color available for each vertex in $H$.
Similarly, all the components of $\left\langle V_2\right\rangle$ can be colored properly. Hence, $\chi(G)\leq \omega(G)+p-1$.
\end{proof}

Even though we are not able to show that the bound given in  Theorem \ref{k1k2kp} is tight, we shall show that the upper bound cannot be made smaller than $\omega+\lceil\frac{p-1}{2}\rceil$ by providing an example that requires $\omega+\lceil\frac{p-1}{2}\rceil$ colors.

 For $p\geq1$, let $X_1= \mathop{\cup}\limits_{i=1}^{p-1}x_i$ ($X_1=\emptyset$ for $p=1$), $X_2= \mathop\cup \limits_{i=p}^{\omega}x_i$, $Y_1= \mathop\cup \limits_{i=1}^{p-1}y_i$ ($Y_1=\emptyset$ for $p=1$), $Y_2= \mathop\cup \limits_{i=p}^{\omega}y_i$ and $Z=\mathop\cup\limits_{i=1}^{p-1}z_i$ ($Z=\emptyset$ for $p=1$). Also, let  $X=X_1 \cup X_2$, $Y=Y_1 \cup Y_2$.
Let $G^*$ be a graph with $V(G^*)=X\cup Y\cup Z$ and edge set $E(G^*)=\{\{x_ix_j\}\cup \{y_iy_j\}\cup \{z_rz_s\}\cup\{x_my_n\}\cup \{y_mz_n\}\cup \{ x_nz_n\}$, where $1\leq i,j,m\leq \omega$, $1\leq r,s,n\leq p-1$, $i\neq j$, $r\neq s$  and $m\neq n$\}. That is,  $X, Y $ and $Z$ form cliques and $[X, Y_1]\cup PM_1$ and $[Y,Z]\cup PM_2$ are complete, where $PM_1= \mathop\cup\limits _{s=1}^{p-1} {x_sy_s}$ and $PM_2=\mathop\cup\limits_{s=1}^{p-1}{y_s z_s}$.  Also, $[X_1,Z]=\{x_s z_s\ |\ 1\leq s\leq p-1\}$, $[X_2,Z]=\emptyset$ and $[X,Y_2]=\emptyset$.
It is easy to observe that $\omega(G^*)=\omega$.
\begin{figure}
  \centering
  \includegraphics[width=10cm]{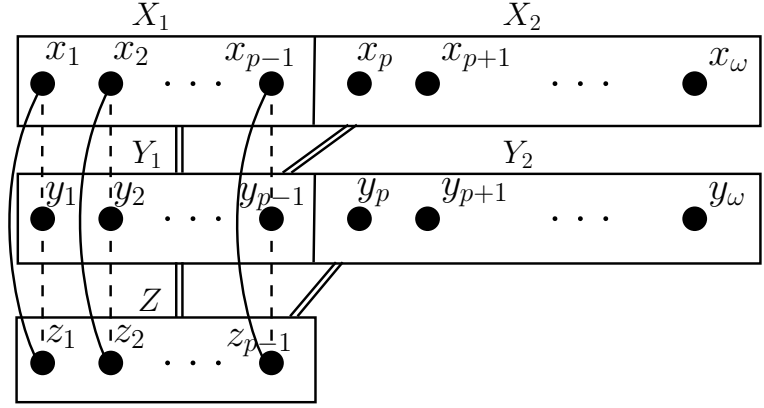}
  \caption{$\{P_3\cup P_2, (K_1\cup K_2)+K_p\}$-free graph}\label{ExK1K2Kp}
\end{figure}

We claim that $G^*$ is a $\{P_3\cup P_2,(K_1 \cup K_2)+K_p\}$-free graph.
Without much difficulty, one can observe that $\alpha(G^*)= 2$.
Also, $\alpha(P_3\cup P_2)=3$.
Hence, $G^*$ is $(P_3 \cup P_2)$-free.
Suppose $G^*$ contains $(K_1\cup K_2)+K_p$ as an induced subgraph.
Since $[Y,Z]\cup PM_2$, $\langle Y\rangle$ and $\langle Z\rangle$ are complete, any three vertices in $Y\cup Z$ will induce a $P_3$ or $K_3$ in $G^*$, and thus the subgraph $K_1\cup K_2$ of $(K_1\cup K_2)+K_p$ cannot be in $\langle Y \cup Z \rangle$.
Hence a vertex $x_i\in X$, for some $i\in \{1,2,\ldots,\omega\}$ should be the part of either $K_1$ or $K_2$ of $(K_1 \cup K_2)+K_p$ in $G^*$.
Similar to the above argument for $Y\cup Z$, the subgraph $K_1\cup K_2$ of $(K_1\cup K_2)+K_p$ cannot be in $\langle X \cup Y_1 \rangle$.
Hence  either a vertex $z_j\in Z$ or $y_k\in Y_2$, for some $j\in \{1,2,\ldots, p-1\}$ and for some $k\in \{p,p+1,\ldots, \omega\}$ should  be a part of  $K_1\cup K_2$ of $(K_1 \cup K_2)+K_p$ in $G^*$.
Note that the common neighbors of the vertices $x_i$ and $z_j$ or $x_i$ and $y_k$ are contained in $Y_1\cup \{x_j\}$ or $Y_1\cup \{z_i\}$ respectively, either of these vertices does not induce a $K_p$ in $G^*$, a contradiction to the fact that $G^*$ contains a $(K_1\cup K_2)+K_p$.
 Hence, $G^*$ is a $\{P_3\cup P_2,(K_1 \cup K_2)+K_p\}$-free graph. Since $\alpha(G^*)=2$, we have $\chi(G^*)\geq\left\lceil\frac{|V(G^*)|}{\alpha(G^*)}\right\rceil=\left\lceil\frac{2\omega(G^*)+p-1}{2}\right\rceil=\omega(G^*)+\left\lceil\frac{p-1}{2}\right\rceil$.

Next,  we obtain a linear $\chi$-binding function for $\{P_3\cup P_2, (K_1\cup K_2)+K_p\}$-free
graphs.
\begin{theorem}\label{k1k2kpallomega}
Let $p$ be an integer greater than $1$. If $G$ is a $\{P_3\cup P_2, (K_1\cup K_2)+K_p\}$-free graph, then $\chi(G)\leq \left\{
\begin{array}{lcl}
\omega(G)+\sum\limits_{j=2}^{p+1} (j-1)(p-j+3)					& \textnormal{for} & 3\leq\omega(G)\leq p+1\\ 
 \omega(G)+7(p-1)+\sum\limits_{j=4}^{p-\left\lfloor\frac{k}{2}\right\rfloor} (j-1)(p-j+3) 	& \textnormal{for} & \omega(G)= (p+2+k), 0\leq k\leq2p-5 \\
 \omega(G)+4p-3  	& \textnormal{for} & \omega(G)= 3p-2\\
\omega(G)+p-1  	& \textnormal{for} & \omega(G)\geq 3p-1.
\end{array}
\right.$
\end{theorem}
\begin{proof}
Let $G$ be a $\{P_3\cup P_2, (K_1\cup K_2)+K_p\}$-free graph with $p\geq2$. For $\omega(G)\geq3p-1$, the bound follows from Theorem \ref{k1k2kp}. By (\ref{emptycij}) of Proposition \ref{k1k2kpprop}, we see that $C_{i,j}=\emptyset$ for all $j\geq p+2$.
We know that $V(G)=V_1\cup V_2$, where $V_1=\mathop\cup\limits_{1\leq \ell\leq\omega}(U_{\ell})$ and $V_2=\mathop\cup\limits_{(i,j)\in L}C_{i,j}$.
Clearly, the vertices of $V_1$ can be colored with $\omega(G)$ colors. Let us find an upper bound for $\chi(\langle V_2\rangle)$.
First, let us consider the case when $3\leq \omega(G)\leq p+1$. For $1\leq i< j\leq p+1$, one can observe that $\omega(\langle C_{i,j}\rangle)\leq \omega(G)-j+2\leq p-j+3$.
Thus one can properly color the vertices of $\left(\mathop\cup\limits_{i=1}^{j-1} C_{i,j}\right)$ with at
most $(j-1)(p-j+3)$ colors and hence $\chi (G) \leq \chi(\langle V_1\rangle)+\chi(\langle V_2\rangle)\leq \omega(G)+\sum\limits_{j=2}^{p+1} (j-1)(p-j+3)$.

Next, let us consider $\omega(G) = (p+2+k)$, where $ 0\leq k\leq2p-4$. By using (\ref{p3p2k1k2kpobs}) of Proposition \ref{p3p2k1k2kpprop},  $\left\langle \left(\mathop\cup \limits_{j=p-\lfloor\frac{k}{2}\rfloor+1}^{p+1} \left(\mathop\cup \limits_{i=1}^{j-1} C_{i,j}\right)\right)\right\rangle$ is a $P_3$-free graph.
By arguments mentioned in Theorem \ref{k1k2kp}, one can color the vertices of $V(G)\backslash \left\{\mathop\cup \limits_{j=2}^{p-\lfloor\frac{k}{2}\rfloor} \left(\mathop\cup \limits_{i=1}^{j-1} C_{i,j}\right)\right\}$ with at most $\omega(G)+p-1$ colors.
For $k=2p-4$, $\omega(G)=3p-2$ and we see that  $\langle V_2\backslash C_{1,2}\rangle$ is $P_3$-free and $\chi(\langle V(G)\backslash C_{1,2}\rangle)\leq\omega(G)+p-1$.
Also,  $\chi(\langle C_{1,2}\rangle)=\omega(\langle C_{1,2}\rangle)\leq 3p-2$ and thus $\chi(G)\leq \omega(G)+4p-3$. Finally, for $0\leq k\leq2p-5$ and $2\leq j\leq {p-\lfloor\frac{k}{2}\rfloor}$, (here $p\geq 3$) let us find the number of colors used for $\left(\mathop\cup\limits_{i=1}^{j-1} C_{i,j}\right)$.
Suppose $\omega (\langle C_{i,j}\rangle)\geq p-j+4$ for some $(i,j)\in L$ such that $1\leq i<j\leq p-\lfloor\frac{k}{2}\rfloor$, by (\ref{omegaCkj}) of
Proposition \ref{p3p2k1k2kpprop}, we see that $\omega (\langle C_{r,j}\rangle)\leq1$ for $r\neq i$ and $1\leq r \leq j-1$.
Also, by (\ref{CijClique}) of Fact \ref{true}, each $C_{i,j}$ is a disjoint union of cliques. In this case, one can properly color the
vertices of $C_{i,j}$  with  $\omega(\langle C_{i,j}\rangle)$ new colors and each $C_{r,j}$ with a new color.
Hence, the number of colors needed for the vertices of $\left(\mathop\cup\limits_{i=1}^{j-1} C_{i,j}\right)$ is at most $\omega (\langle C_{i,j}\rangle)+(j-2)\leq(\omega(G)-j+2)+(j-2)=\omega(G)$. Suppose $\omega(\langle C_{i,j}\rangle)\leq (p-j+3)$ for every $(i,j)\in L$ such that
$1\leq i<j\leq p-\lfloor\frac{k}{2}\rfloor$, then one can properly color the vertices of $\left(\mathop\cup\limits_{i=1}^{j-1} C_{i,j}\right)$ with at
most $(j-1)(p-j+3)$ colors.
Now, for $4\leq j\leq p-\left\lfloor\frac{k}{2}\right\rfloor$ let us compare the values of $(j-1)(p-j+3)$ and $\omega (G)$ . \\
$(j-1)(p-j+3)-\omega(G)$ $\geq(j-1)p-(j-1)(j-3)-(3p-3)$\\
\phantom{$(j-1)(p-j+3)-\omega(G)$}$=(j-4)p-(j-1)(j-4)-(j-1)+3$\\
\phantom{$(j-1)(p-j+3)-\omega(G)$}$=(j-4)(p-j)+j-4-j+4$\\
\phantom{$(j-1)(p-j+3)-\omega(G)$}$=(j-4)(p-j)\geq 0$.\\
Thus $(j-1)(p-j+3)\geq\omega(G)$ and hence the vertices of $\left(\mathop\cup \limits_{j=4}^{p-\lfloor\frac{k}{2}\rfloor} \left(\mathop\cup \limits_{i=1}^{j-1} C_{i,j}\right)\right)$ can be properly colored with $\sum\limits_{j=4}^{p-\lfloor\frac{k}{2}\rfloor} (j-1)(p-j+3)$ colors.
For $p\geq 3$ and $j=2,3$, $(j-1)(p-j+3)\leq 3p-3$.
Therefore the vertices of $C_{1,2}$ and $(C_{1,3}\cup C_{2,3})$ can be properly colored with at most $(3p-3)$ colors each.
Hence,  $\chi (G) \leq (\omega(G)+p-1)+2(3p-3)+\sum\limits_{j=4}^{p-\lfloor\frac{k}{2}\rfloor} (j-1)(p-j+3)=\omega(G)+7(p-1)+\sum\limits_{j=4}^{p-\lfloor\frac{k}{2}\rfloor} (j-1)(p-j+3)$.
\end{proof}

The $\chi$-binding function obtained in Theorem \ref{k1k2kpallomega} is not optimal.
For instance, if $p=2$, $(K_1+K_2)+K_p\cong HVN$.
For $\omega(G)\geq 4$, we shall show that any $\{P_3\cup P_2,HVN\}$-free graph is $(\omega +1)$-colorable.
This can be seen in Theorem \ref{hvn}.

\begin{theorem}\label{hvn}
If $G$ is a $\{P_3\cup P_2, HVN\}$-free graph with $\omega(G)\geq4$, then $\chi(G)\leq\omega(G)+1$.
\end{theorem}
\begin{proof}
For $\omega(G)\geq 5$, the bound follows from Theorem \ref{k1k2kp}.
The only remaining case is $\omega(G)=4$. By Proposition \ref{k1k2kpprop}, we see that $V(G)=V_1\cup V_2$ where
$V_1=\mathop\cup\limits_{1\leq k\leq 4} U_k$ and $V_2=C_{1,2}\cup C_{1,3}\cup C_{2,3}$ and each vertex in $V_2$ is
adjacent to the vertices of at most one $U_{\ell}$, $\ell\in\{1,2,3,4\}$. By (\ref{CijClique}) of Fact \ref{true},
each $\langle C_{i,j}\rangle$ is $P_3$-free, where  $1\leq i<j\leq 3$.
We shall show that $G$ is $\omega(G)+1=5$ colorable.
Let $\{1,2,3,4,5\}$ be the set of colors. For $1\leq k\leq 4$, let us assign the color $k$ to the vertices of $U_k$.
In order to color the vertices of $V_2$, let us break the proof into two cases depending upon the presence of edges in $\langle C_{i,j}\rangle$, $1\leq i<j\leq3$.

\noindent\textbf{Case 1} $\langle C_{1,3}\rangle$ or $\langle C_{2,3}\rangle$ contains an edge.

First, let us assume $\langle C_{2,3}\rangle$ contains an edge, say $ab$.
We claim that both $C_{1,2}$ and $C_{1,3}$ are independent sets.
Suppose, there exists an edge $cd \in \langle C_{1,j}\rangle$ for $ j\in\{2,3\}$.
First, we show that $[\{a,b\},\{c,d\}]$ is neither empty nor complete.
If $[\{a,b\},\{c,d\}]=\emptyset$, then $\langle\{a,v_1,v_j,c,d\}\rangle\cong P_3\cup P_2$, a contradiction.
If $[\{a,b\},\{c,d\}]$ is complete, then $\langle\{a,b,c,d,v_1\}\rangle\cong HVN$, a contradiction.
Hence, there exists an induced $P_3\in \langle\{a,b,c,d\}\rangle$, say $P$. If $cd\in \langle C_{1,3}\rangle$, then
$\langle\{V(P),v_3,v_4\}\rangle\cong P_3\cup P_2$, a contradiction. Thus, $C_{1,3}$ is an independent set.
Let us consider $cd\in \langle C_{1,2}\rangle$.
If $|N_{\{v_3,v_4\}}(V(P))|\leq 1$, then $N_{\{v_3,v_4\}}(V(P))\subseteq\{v_\ell\}$, $\ell\in \{3,4\}$, and hence   $\langle\{V(P),v_2,$ $\{v_3,v_4\}\backslash\{v_{\ell}\}\}\rangle\cong P_3\cup P_2$, a contradiction.
Suppose $|N_{\{v_3,v_4\}}(V(P))|=2$, by (\ref{p-1neighbors}) of Proposition \ref{k1k2kpprop}, $|N_A(c)|=|N_A(d)|=1$
and both $c,d\in V(P)$.
Without loss of generality, let us assume that $cv_3,dv_4\in E(G)$ and $V(P)=\{a,c,d\}$ such that $ad\notin E(G)$. Now, we have $\langle\{v_3,v_4,d,a,b\}\rangle\cong P_3\cup P_2$ when $bd\notin E(G)$ and
$\langle\{a,b,d, v_2,v_3\}\rangle\cong P_3\cup P_2$ when $bd\in E(G)$, a contradiction.
Therefore $\langle C_{1,2}\rangle$ is also an independent set.
Similarly when $\langle C_{1,3}\rangle$ contains an edge, we can show that both $C_{1,2}$ and $C_{2,3}$ are independent sets.

Let $\{i,k\}=\{1,2\}$. If $\langle C_{i,3}\rangle$ contains an edge, then we can establish a $5$-coloring for $G$
by assigning the colors from $\{i,3,4\}$, $5$ and $k$ to the vertices of $C_{i,3}$, $C_{1,2}$ and $C_{k,3}$ respectively.
Since $\langle C_{i,3}\rangle$ is $P_3$-free, $\omega(\langle C_{i,3}\rangle)\leq 3$, $C_{1,2}$ and $C_{k,3}$
are independent sets, we observe that this is a proper $5$-coloring for $G$.

\noindent\textbf{Case 2} $\langle C_{1,2}\rangle$ contains an edge.

Here, by Case 1, $C_{1,3}$ and $C_{2,3}$ are independent sets.
Let $S\subseteq  C_{1,2}$ be the vertices of a maximum clique in  $\langle C_{1,2}\rangle$. Let us
break this case into three subcases depending upon the size of the maximum clique $S$.

\noindent\textbf{Case 2.1} $|S|= 2$

Each vertex in $C_{1,2}$ is adjacent to at most one color in $\{1,2,3,4\}$.
Let us assign the available color from $\{1,2,5\}$ to the vertices of  $C_{1,2}$ and assign the color $3$
and $4$ to the vertices of $C_{1,3}$ and $C_{2,3}$ respectively.
Since $\langle C_{1,2}\rangle$ is $P_3$-free, $C_{1,3}$ and $C_{2,3}$ are independent sets, this will be a proper coloring for $G$.

\noindent\textbf{Case 2.2} $|S|=3$

If $I_k=\emptyset$ for every $k\in\{1,2,3,4\}$, then a coloring similar to that given in case 2.1 will give us a proper 5-coloring for $G$.
So, let us assume that $I_k\neq \emptyset$ for some $k\in\{1,2,3, 4\}$ and let $u_k\in I_k$ and $S=\{a,b,c\}$.
We show that $|N_S(\{v_k,u_k\})|\geq 2$.
Suppose $|N_S(\{v_k,u_k\})|\leq 1$, without loss of generality, let $b,c\notin N_S(\{v_k,u_k\})$.
Let $\ell\in \{1,2\}\backslash \{k\}$.
Now, we have $\langle\{u_k,v_{\ell},v_k,b,c\}\rangle\cong P_3\cup P_2$, a contradiction.
Next, we show that $I_q=\emptyset$, for every $q\in\{1,2,3,4\}$ and $q\neq k$.
Suppose $u_q\in I_q$, by using the above arguments we have $|N_S(\{v_q, u_q\})|\geq 2$ which implies that there exists a vertex in $S$ which has a neighbor in both $\{v_k, u_k\}\subseteq U_k$
and $\{v_q, u_q\}\subseteq U_q$, a contradiction.

Now, let us establish a proper $5$-coloring for $G$. 
If $|N_S(\{v_k, u_k\})|=2 \text{\ or\ }3$ and the vertices of $S$ are adjacent to only vertices of color $k$,
then color the vertices of $S$ with the colors $\{1,2,3,4\}\backslash \{k\}$.
If $|N_S(\{v_k, u_k\})|=2$ and the vertices of $S$ are adjacent only to the vertices of colors $k$ and $q$, then
assign the color $q$ to a vertex of $N_S(\{v_k, u_k\})$ and the colors $\{1,2,3,4\}\backslash\{k,q\}$ for the remaining vertices of $S$.
Similarly, we can color all the components of $\langle C_{1,2}\rangle$ by using the colors $\{1,2,3,4\}\backslash \{k\}$.
Now, if $k=1$, then assign the colors $k$ and $5$ to the vertices of $C_{1,3}$ and $C_{2,3}$ respectively.
If $k\neq 1$, then assign the colors $k$ and $5$ to the vertices of $C_{2,3}$ and $C_{1,3}$ respectively.
Clearly, this is a proper coloring for $G$.

\noindent\textbf{Case 2.3} $|S|=4$

Let $S=\{a,b,c,d\}$. We know that $C_{1,3}$ and $C_{2,3}$ are independent. In addition, we shall show that $(C_{1,3}\cup C_{2,3})$ is also an independent set.
On the contrary, let us assume that there exists an edge $ef\in E(\langle C_{1,3}\cup C_{2,3}\rangle)$, say  $e\in C_{1,3}$ and $f\in C_{2,3}$.
Since $\omega(G)=4$, we have $|N_S(e)|\leq3$ and $|N_S(f)|\leq3$. 
Suppose $|N_S(e)|\leq 2$, then there exists two vertices of $S$ which are not in  $N_S(e)$, say $c,d$ and hence 
$\langle\{e,v_2,v_1,c,d\}\rangle\cong P_3\cup P_2$, a contradiction.
Thus, the only possibility is $|N_S(e)|=3$.
Similarly, we get $|N_S(f)|=3$. Thus, $|N_S(e)\cap N_S(f)|\geq 2$.
Without loss of generality, let us assume $a,b\in (N_S(e)\cap N_S(f))$.
Since $|N_A(d)|\leq 1$, either $dv_3\notin E(G)$ or $dv_4\notin E(G)$. Without loss of generality,
let us assume $dv_3\notin E(G)$ and $N_S(e)=\{a,b,c\}$.
If $xv_3\notin E(G)$ for some $x\in \{a,b\}$, then $\langle\{e,x,d,v_1,v_3\}\rangle\cong P_3\cup P_2$, a contradiction.
Therefore $[\{a,b\},v_3]$ is complete, but even here we get $\langle\{a,b,e,f,v_3\}\rangle\cong HVN$, a contradiction.
Hence, $(C_{1,3}\cup C_{2,3})$ is an independent set.
Next, we show that $N_S(x)=N_S(y)$ for any $x,y\in (C_{1,3}\cup C_{2,3})$.
Suppose, $N_S(x)\neq N_S(y)$ for some $x,y\in (C_{1,3}\cup C_{2,3})$, by similar arguments we see that $|N_S(x)|=|N_S(y)|=3$.
Without loss of generality, let us assume that $N_S(x)=\{a,b,c\}$ and $N_S(y)=\{a,b,d\}$. Here
$\langle\{a,b,c,x,y\}\rangle\cong HVN$, a contradiction.

Let us establish a $5$-coloring for $G$.
We first consider the case when $I_k=\emptyset$, for all $1\leq k\leq 4$.
Since $\omega(G)=4$, there exists a vertex $x\in S$ such that $xv_3\notin E(G)$.
Let us assign the color $3$ to the vertex $x$ and  color the remaining vertices of $S$ with the colors $\{1,2,5\}$.
Similarly, we can color all the components of $C_{1,2}$ with the colors $\{1,2,3,5\}$. Finally assign the color $4$ to the vertices of $(C_{1,3}\cup C_{2,3})$. This yields a 5-coloring for $G$.

Let us next consider $I_k\neq \emptyset$, for some $k\in\{1,2,3,4\}$ and let $u_k\in I_k$. By arguments similar to that given in case $2.2$, we can show that $|N_S(\{v_k, u_k\})|\geq3$ and as a consequence $I_q=\emptyset$, for all $1\leq q\neq k\leq 4$ .
Let $k\in \{3,4\}$. If $|N_S(\{v_k, u_k\})|=3 \text{\ or\ }4$ and the vertices of $S$ are adjacent only  to
 vertices with color $k$, then color the vertices of $S$ with the colors $\{1,2,3,4,5\}\backslash \{k\}$.
If $|N_S(\{v_k, u_k\})|=3$ and the vertices of $S$ are adjacent only to vertices with colors $k$ and $q$, then assign the color $q$ to a vertex of $N_S(\{v_k, u_k\})$ and color the remaining vertices of $S$ with the colors $\{1,2,3,4,5\}\backslash \{k,q\}$.
Similarly, we can color all the components of $\langle C_{1,2}\rangle$ with the colors $\{1,2,3,4,5\}\backslash\{k\}$. Finally, assign the color $k$
to the vertices of $(C_{1,3}\cup C_{2,3})$. One can observe that, this is a 5-coloring of $G$.
Next let us consider $k\in \{1,2\}$. Clearly, $N_S(v_k)=\emptyset$ and thereby $|N_S(\{u_k\})|\geq3$, for any $u_k\in I_k$. Since $\omega(G)=4$, $|N_S(\{u_k\})|=3$ for any $u_k\in I_k$. 
Now, we show that $|N_S(I_k)|= 3$.
Suppose $|N_S(I_k)|=4$, there exists at least two vertices $x,y\in I_k$ such that $N_S(x)\neq N_S(y)$.
Without loss of generality, let us assume $N_S(x)=\{a,b,c\}$ and $N_S(y)=\{a,b,d\}$, then
$\langle\{a,b,c,x,y\}\rangle\cong HVN$, a contradiction.
Hence, $|N_S(I_k)|=3$ and there exists a vertex say $d\in S$ such that $d\notin N_S(\{v_k,I_k\})$.
Let us assign the color $k$ to the vertex $d$ and color the remaining vertices of $S$ with
the colors $\{1,2,3,4\}\backslash\{k\}$.
Similarly, we can color all the components of $\langle C_{1,2}\rangle$.
Finally, assign the color $5$ to the vertices of $(C_{1,3}\cup C_{2,3})$.
Clearly, this is a proper $5$-coloring for $G$.
\end{proof}

For $\omega(G)\geq4$, the bound given in Theorem \ref{hvn} is tight. Let $G^*$ be the graph as shown in Figure \ref{ExK1K2Kp} with $\omega(G^*)\geq4$ and $p=2$.
We have already observed that $G^*$ is $\{P_3\cup P_2,HVN\}$-free. Clearly, $|V(G^*)|=2\omega(G^*)+1$ and $\alpha (G^*)=2$. Hence, $\chi(G^*)\geq\left\lceil \frac{V(G^*)}{\alpha (G^*)}\right\rceil = \left\lceil\frac{2\omega(G^*)+1}{2}\right\rceil=\omega(G^*)+1$.


\section{{$\{P_3\cup P_2, 2K_1+K_p\}$}-free graphs}

Any $\{P_3\cup P_2, 2K_1+K_p\}$-free graphs are also a $\{P_3\cup P_2, (K_1\cup K_2)+K_p\}$-free graphs.
Hence the properties established for $\{P_3\cup P_2, (K_1\cup K_2)+K_p\}$-free graphs are also true for $\{P_3\cup P_2, 2K_1+K_p\}$-free graphs. 
Note that when $p=0$ or $1$, $\{P_3\cup P_2, 2K_1+K_p\}$ is $P_3$-free and hence perfect.
Therefore we shall assume that $p\geq 2$.
Let us begin Section 4, by obtaining a structural characterization and the chromatic number of  $\{P_3\cup P_2,2K_1+K_p\}$-free graph $G$, when $\omega(G)\geq 3p-1$.
\begin{theorem}\label{2k1+kp3p-1}
Let $G$ be a $\{P_3\cup P_2, 2K_1+K_p\}$-free graph with $p\geq 2$ and $V(G)=V_1\cup V_2$. If $\omega(G)\geq 3p-1$, then $\langle V_1\rangle$ is complete, $\langle V_2\rangle$ is $P_3$-free and $G$ is perfect.
\end{theorem}

\begin{proof}
Let $G$ be a $\{P_3\cup P_2, 2K_1+K_p\}$-free graph with $\omega(G)\geq 3p-1$, $p\geq 2$.
First, we shall observe that $I_{k}=\emptyset$  for every  $k\in \{1,2,\ldots,\omega(G)\}$.
On the contrary, let $u_k\in I_k$ for some $k\in \{1,2,\ldots,\omega\}$.
Since $p\geq2$ and $\omega(G)\geq 3p-1\geq p+1$, the vertices $u_k,v_k$ and any $p$ vertices of $A\backslash\{v_k\}$ will induce a $(2K_1+K_p)$, a contradiction.
Thus $V_1=A$ and hence complete. 
Since any $\{P_3\cup P_2, 2K_1+K_p\}$-free graph is also a $\{P_3\cup P_2, (K_1\cup K_2)+K_p\}$-free graph, by Theorem \ref{k1k2kp}, we see that  $\langle V_2\rangle$ is $P_3$-free.
Finally, to show that $G$ is perfect, by using Theorem \ref{SPGT}, it is enough to show that for any $r\geq 2$, $C_{2r+1}$ and $\overline C_{2r+1}$ are not induced subgraphs of $G$.
Since $G$ is $(P_3\cup P_2)$-free, $C_{2q+1}\notin G$, $q\geq  3$.
Note that $\overline{C}_5\cong C_5$ and $\omega(\overline C_{2r+1})=r$. Suppose there exists $\overline C_{2r+1}\in G$, $r\geq 2$,  then $\overline C_{2r+1}$ can have at most $r$ vertices in $A$ and hence must have at least $(r+1)$ vertices in $V_2$ which will induced a $P_3$, a contradiction.
\end{proof}

As a consequences of Theorem \ref{k1k2kpallomega} and Theorem \ref{2k1+kp3p-1}, we obtain Proposition \ref{2k1kpomegageq2p} and Corollary \ref{2k1kpallomega}.
\begin{proposition}\label{2k1kpomegageq2p}
Let $G$ be a $\{P_3\cup P_2, 2K_1+K_p\}$-free graph with $p\geq 2$. If $\omega(\langle C_{1,2}\rangle)\geq 2p$, then $\chi(G)=\omega(G)$.
\end{proposition}
\begin{proof}
Let $G$ be a $\{P_3\cup P_2, 2K_1+K_p\}$-free graph with $p\geq 2$ and $\omega(\langle C_{1,2}\rangle)\geq 2p$.
Since $\omega(G)\geq\omega(\langle C_{1,2}\rangle)\geq 2p\geq p+1$, as in Theorem \ref{2k1+kp3p-1}, we see that $V(G)=V_1\cup V_2$, where $V_1=A$ and
$V_2=\mathop\cup\limits_{(i,j)\in L}C_{i,j}$. Now if $\omega(G)\geq 3p-1$, by Theorem \ref{2k1+kp3p-1}, we see that $\langle V_2\rangle$ is $P_3$-free and $\chi(G)=\omega(G)$.
Therefore it is enough to prove that $\chi(G)=\omega(G)$ when $2p\leq \omega(G)\leq  3p-2$.
Let us consider $S\subseteq C_{1,2}$ such that $\langle S\rangle$ is a maximum clique in $\langle C_{1,2}\rangle$.
For any vertex $a\in V_2\backslash C_{1,2}$, we shall show that $N_S(a)=S$.
Suppose there exists a vertex $a\in V_2\backslash C_{1,2}$ such that $N_S(a)\neq S$, then there exists a vertex $u\in S$ such that $ua\notin E(G)$.
If $|N_S(a)|\geq p$, then there exist $p$ vertices of $N_S(a)$ say $S'$ such that $\langle\{a,u\}\cup S'\rangle\cong 2K_1+K_p$, a contradiction.
Thus $|N_S(a)|\leq p-1$.
By (\ref{p-1neighbors}) of Proposition \ref{k1k2kpprop}, we have $|N_A(a)|\leq p-1$ and thus $|A\backslash N_A(a)|\geq p+1$ as $\omega(G)\geq 2p$.
Since $\omega(G)\leq 3p-2$ and $|S|\geq2p$, there exists a vertex $v_r\in A\backslash N_A(a)$ such that $N_S(v_r)\neq S$ (otherwise $\omega(\langle S\cup \{A\backslash N_A(a)\}\rangle)\geq 3p+1$, a contradiction).
For the same reason as mentioned for $|N_S(a)|\leq p-1$, we see that $|N_S(v_r)|\leq p-1$.
Therefore $|N_S(\{a,v_r\})|\leq 2(p-1)$ and as a consequence there exist  at least two vertices $x,y\in S$ such that $[\{a,v_r\},\{x,y\}]=\emptyset$.
Since $a\in V_2\backslash C_{1,2}$, $av_t\in E(G)$ for some $t\in \{1,2\}$, and thereby $\langle\{a,v_t,v_r,x,y\}\rangle\cong P_3\cup P_2$, a contradiction.
Thus, $N_S(a)=S$  for any vertex $a\in V_2\backslash C_{1,2}$.

Since $G$ is $(2K_1+K_p)$-free and $|S|\geq2p$, we see that  $\langle V_2\backslash C_{1,2}\rangle$ is $2K_1$-free. Hence $\langle S\cup \{V_2\backslash C_{1,2}\}\rangle$ is a clique and thus it is $P_3$-free.
By arguments mentioned in Theorem \ref{2k1+kp3p-1}, we get that $\langle A\cup S\cup \{V_2\backslash C_{1,2}\}\rangle$ is perfect and these vertices can be colored with $\omega(G)$ colors $\{1,2,\ldots,\omega(G)\}$ such that the color $i$ is assigned to vertex $v_i$, $1\leq i\leq \omega(G)$.
Now, we show that each of the color given to the vertices of $S$ is available for every vertex in $C_{1,2}\backslash S$.
Suppose there exists a color $q$ assigned to a vertex in $S$ which is unavailable to some vertex $b\in\langle C_{1,2}\backslash S\rangle$ then $N_S(v_q)\neq S$ and $bv_q\in E(G)$.
Again for the same reason as mentioned for $|N_S(a)|\leq p-1$, we have $|N_S(v_q)|\leq p-1$ and there exist two vertices $x,y\in S$ such that $x,y\notin N_S(v_q)$.
Now, we get $\langle\{ b,v_q,v_1,x,y\}\rangle\cong P_3\cup P_2$, a contradiction.
Since $S$ induce a maximum clique in $\langle C_{1,2}\rangle$, the vertices in $C_{1,2}\backslash S$ can be colored using the colors assigned to the vertices of $S$.
Hence $\chi(G)=\omega(G)$.
\end{proof}
\begin{corollary}\label{2k1kpallomega}
Let $G$ be a $\{P_3\cup P_2, 2K_1+K_p\}$-free graph with $p\geq 2$, then\\ $\chi(G)\leq \left\{
\begin{array}{lcl}
\omega(G)+\sum\limits_{j=2}^{p+1} (j-1)(p-j+3)					& \textnormal{for} & 3\leq\omega(G)\leq p+1\\  \omega(G)+2p-1+\sum\limits_{j=3}^{p-\left\lfloor\frac{k}{2}\right\rfloor} (j-1)(p-j+3) 	& \textnormal{for} & \omega(G)= (p+2+k), 0\leq k\leq2p-5 \\
 \omega(G)+2p-1  	& \textnormal{for} & \omega(G)= 3p-2\\
\omega(G)	& \textnormal{for} & \omega(G)\geq 3p-1.
\end{array}
\right.$
\end{corollary}
\begin{proof}
Let $G$ be $\{P_3\cup P_2, 2K_1+K_p\}$-free graph with $p\geq 2$.
When $3\leq\omega(G)\leq p+1$ or when  $\omega(G)\geq 3p-1$, the bound follows from Theorem \ref{k1k2kpallomega} and Theorem \ref{2k1+kp3p-1} respectively.
Hence, let us consider $\omega(G)= (p+2+k), 0\leq k\leq2p-4$.
As in Theorem \ref{2k1+kp3p-1}, we see that $V(G)=V_1\cup V_2$,
where $V_1=A$ and $V_2=\mathop\cup\limits_{(i,j)\in L}C_{i,j}$.
Also, by Proposition \ref{2k1kpomegageq2p}, $\chi(G)=\omega(G)$ when $\omega(\langle C_{1,2}\rangle)\geq2p$.
Hence, we shall consider $\omega(\langle C_{1,2}\rangle)\leq 2p-1$.
By (\ref{p3p2k1k2kpobs}) of Proposition \ref{p3p2k1k2kpprop}, we see that $\left\langle \left(\mathop\cup \limits_{j=p-\lfloor\frac{k}{2}\rfloor+1}^{p+1} \left(\mathop\cup \limits_{i=1}^{j-1} C_{i,j}\right)\right)\right\rangle$ is a $P_3$-free graph.
As done in Theorem \ref{2k1+kp3p-1}, we can observe that $\left\langle A\cup \left(\mathop\cup \limits_{j=p-\lfloor\frac{k}{2}\rfloor+1}^{p+1} \left(\mathop\cup \limits_{i=1}^{j-1} C_{i,j}\right)\right)\right\rangle$ is perfect and thus we can color the vertices of $V(G)\backslash \left\{\mathop\cup \limits_{j=2}^{p-\lfloor\frac{k}{2}\rfloor} \left(\mathop\cup \limits_{i=1}^{j-1} C_{i,j}\right)\right\}$ with $\omega(G)$ colors.
In particular, for $\omega(G)=3p-2$, $\langle V_2\backslash C_{1,2}\rangle$ is $\omega(G)$-colorable and thus $\chi(G)\leq \omega(G)+2p-1$.
Finally, the cases remaining are $\omega(G) = (p+2+k)\geq4, 0\leq k\leq2p-5$ when $\omega(\langle C_{1,2}\rangle)\leq 2p-1$.
For the same reason as given in  Theorem \ref{k1k2kpallomega}, the vertices of $\left(\mathop\cup \limits_{j=4}^{p-\lfloor\frac{k}{2}\rfloor} \left(\mathop\cup \limits_{i=1}^{j-1} C_{i,j}\right)\right)$ can be properly colored with $\sum\limits_{j=4}^{p-\lfloor\frac{k}{2}\rfloor} (j-1)(p-j+3)$ colors.
Hence the vertices of $V(G)\backslash (C_{1,2}\cup C_{1,3}\cup C_{2,3})$ can be properly colored with at most $\omega(G)+\sum\limits_{j=4}^{p-\lfloor\frac{k}{2}\rfloor} (j-1)(p-j+3)$ colors and therefore $\chi(G\backslash (C_{1,3}\cup C_{2,3}))\leq \omega(G)+2p-1+\sum\limits_{j=4}^{p-\lfloor\frac{k}{2}\rfloor}(j-1)(p-j+3)$.
All that remains is to find the number of new colors required to color the vertices of $C_{1,3}\cup C_{2,3}$.
Let $\{i,k\}=\{1,2\}$.
First let us observe that if $\omega(\langle C_{i,3}\rangle)\geq p+1$, then $C_{k,3}=\emptyset$.
On the contrary let $\omega(\langle C_{i,3}\rangle)\geq p+1$ and $C_{k,3}\neq\emptyset$.
Let $a\in C_{k,3}$ and 
$S\subseteq C_{i,3}$ be vertices which induce a clique of size $p+1$. 
If $|N_S(a)|\geq p$, then there exist $p$ vertices of $N_S(a)$ say $S'$ such that  $\langle \{a,v_k\}\cup S'\rangle\cong 2K_1+K_p$, a contradiction.
If $|N_S(a)|\leq p-1$, then there exist two vertices $x,y\in S$ but not in $N_S(a)$ and clearly $\langle\{ a,v_i,v_3,x,y\}\rangle\cong P_3\cup P_2$, a contradiction. 
Hence if $\omega(\langle C_{i,3}\rangle)\geq p+1$, then $C_{k,3}=\emptyset$.
Further we see that if $\omega(\langle C_{i,3}\rangle)\geq 2p$, then by similar arguments as in Proposition \ref{2k1kpomegageq2p}, one can color the vertices of $V(G)\backslash C_{1,2}$ with $\omega(G)$ colors and hence no new color is necessary for coloring the vertices in $C_{1,3}\cup C_{2,3}$.
If $p+1\leq \omega(\langle C_{i,3}\rangle)\leq 2p-1$, then $C_{1,3}\cup C_{2,3}$ can be colored with at most $2p-1$ new colors.
Suppose $\omega(\langle C_{1,3}\rangle)\leq p$ and $\omega(\langle C_{2,3}\rangle)\leq p$, the vertices of $C_{1,3}\cup C_{2,3}$ can be colored with at most $2p$ new colors.
Thus in any case, the vertices of $C_{1,3}\cup C_{2,3}$ can be colored with at most $2p$ new colors.
Also, $(j-1)(p-j+3)=2p$ when $j=3$.
Hence, $\chi(G)\leq \omega(G)+2p-1+\sum\limits_{j=3}^{p-\lfloor\frac{k}{2}\rfloor} (j-1)(p-j+3)$.
\end{proof}

When $p=2$, we see that $2K_1+K_p\cong diamond$ and hence by Theorem \ref{2k1+kp3p-1} we see that $\{P_3\cup P_2, diamond\}$-free graphs are perfect for $\omega(G)\geq 5$. This was shown by A. P. Bharathi and S. A. Choudum in \cite{bharathi2018colouring}.

\begin{theorem}(\cite{bharathi2018colouring}) \label{p3p2diamond}
If $G$ is a $\{P_3\cup P_2, diamond\}$-free graph then $\chi(G)\leq \left\{
\begin{array}{lcl}
4	& \textnormal{for} & \omega(G)= 2\\
6 & \textnormal{for} & \omega(G)= 3\\
5 & \textnormal{for} & \omega(G)= 4\\
\end{array}
\right.$
and $G$ is perfect if $\omega(G)\geq 5$.
\end{theorem}

We further improve the bound given in Theorem \ref{p3p2diamond}  by obtaining a $\omega(G)$-coloring when $\omega(G)=4$.
This can be seen in Theorem \ref{diamond}.
\begin{theorem}\label{diamond}
If $G$ is a $\{P_3\cup P_2, diamond\}$-free graph then $\chi(G)\leq \left\{
\begin{array}{lcl}
4	& \textnormal{for} & \omega(G)= 2\\
6 & \textnormal{for} & \omega(G)= 3\\
4 & \textnormal{for} & \omega(G)=4.
\end{array}
\right.$
and $G$ is perfect if $\omega(G)\geq 5$.
\end{theorem}
\begin{proof}
For $\omega(G)=2$, $\omega(G)=3$ and $\omega(G)\geq 5$, the results follows from Theorem \ref{p3p2diamond} and \ref{2k1+kp3p-1}.
The only remaining case is $\omega(G)=4$.
Let us establish a 4-coloring for $G$ using the colors $\{1,2,3,4\}$.
By Proposition \ref{k1k2kpprop} and Theorem \ref{2k1+kp3p-1}, we see that $V(G)=V_1\cup V_2$ such that
$V_1=A$, $V_2=C_{1,2}\cup C_{1,3}\cup C_{2,3}$ and each vertex in $V_2$ has
at most one neighbor in $A$. By (\ref{joinvertex}) of Fact \ref{true},
each $\langle C_{i,j}\rangle$ is $P_3$-free, $1\leq i<j\leq 3$.
For $1\leq k\leq 4$, let us assign the color $k$ to the vertex $v_k$.
Let $S$ be the vertices of a maximum clique in $\langle C_{1,2}\rangle$.
Now let us color $V_2$ by considering the possibility of $|S|$.

\textbf{Case 1} $|S|=4$  

It follows from Proposition \ref{2k1kpomegageq2p} that $G$ is $4$-colorable.

\textbf{Case 2} $|S|=3$

Here we first observe that for any vertex $x\in C_{1,3}\cup C_{2,3}$, $[x,S]$ is complete.
Suppose $[x,S]$ is not complete, then there exists a vertex $a\in S$ such that $xa\notin E(G)$.
Let $b\in S \backslash \{a\}$. If $xb\notin E(G)$, then $\langle\{\{x,v_1,v_2\},a,b\}\rangle \cong P_3 \cup P_2$.
Else, $\langle S\cup \{x\}\rangle\cong diamond$, a contradiction.
Thus $[x,S]$ is complete.
As a consequence we see that $|C_{1,3}\cup C_{2,3}|\leq 1$ (otherwise $\langle S\cup C_{1,3}\cup C_{2,3} \rangle$) will contain an induced $diamond$ or a $K_5$, a contradiction).
Now, color the vertices of $S$ properly with the colors $\{1,2,q\}$, where  $q\in \{3,4\}$ and $C_{1,3}\cup C_{2,3}$ with $\{3,4\}\backslash \{q\}$.
Next, let us show that the colors $\{1,2,q\}$ are available for every vertex in $C_{1,2}\backslash S$.
Suppose there exist a vertex $u \in C_{1,2}  \backslash S $ such that u cannot be colored with $q$ then $uv_q\in E(G)$.
Since the color $q$ is given to some vertex in $S$ and $v_q$ being adjacent to two vertices in $S$ will yield a $diamond$, $|[v_q,S]|\leq 1$.
Hence $[v_q,\{a,b\}]=\emptyset$ for some $a,b\in S$ and therefore $\langle u,v_q,v_1,a,b\rangle\cong P_3\cup P_2$, a contradiction.
Thus we get a proper $4$-coloring for $G$.

\textbf{Case 3} $|S|\leq 2$

Here if $\omega (\langle C_{1,3}\rangle)\leq 1$ and $\omega (\langle C_{2,3}\rangle)\leq 1$, then we can properly color the vertices of $C_{1,2},C_{1,3}$ and $C_{2,3}$ with the colors $\{1,2\},3$ and $4$ respectively. 
Suppose $\langle C_{1,3}\rangle$ or $\langle C_{2,3}\rangle$ contains an edge. 
Without loss of generality, let $ab \in E(\langle C_{1,3}\rangle)$. 
As in Case 1 of Theorem \ref{hvn} we can see that $C_{1,2}$ and $C_{2,3}$ are independent sets.
In addition, we can show that $C_{2,3}=\emptyset$.
Suppose $C_{2,3}\neq \emptyset$, let $x \in C_{2,3}$.
If $\langle \{x,a,b\}\rangle\cong K_3$, then $\langle \{x,a,b,v_2\}\rangle\cong diamond$. 
If $\langle\{x,a,b\}\rangle\cong P_3$, then $\langle\{x,a,b,v_3,v_4\}\rangle\cong P_3 \cup P_2$. 
If $[x,\{a,b\}]=\emptyset$, then $\langle\{x,v_1,v_3,a,b\}\rangle\cong P_3 \cup P_2$.
Hence we get a contradiction in all the cases. 
Therefore, $C_{2,3}=\emptyset$ and $|S|\leq1$. 
Since $\omega(G)=4$, we see that $\omega(\langle C_{1,3}\rangle)\leq3$ and thereby we can properly color the vertices of $C_{1,2}$ and $C_{1,3}$ with the colors $2$ and $\{1,3,4\}$ respectively.
Thus we have a proper $4$-coloring for $G$.
\end{proof}

Let $G$ be the complement of the Schl\"{a}fli graph (see, \textcolor{blue}{www.distanceregular.org/graphs/complement-schlafli.html}). Note that $G$ is a $10$-regular graph with $27$ vertices such that any two adjacent vertices have exactly one common neighbor (see, \textcolor{blue}{http://en.wikipedia.org/wiki/Locally\_linear\_graph}) and hence $G$ is $diamond$-free. Next, we shall observe that $G$ is $(P_3\cup P_2)$-free.
Let $uv\cong P_2$ in $G$. The vertices $u$ and $v$ have exactly one common neighbor, say $w$. The ends of the edges $uv, vw$ and $wu$ does not have a common neighbor other than $w,u$ and $v$ respectively. Since $G$ is $10$-regular, each of $u,v$ and $w$ should have exactly $8$ distinct neighbors. Clearly, $|V(G)|-|N({u,v})|=27-(3+8+8)=8$. Thus $V(G)\backslash\{N({u,v})\}=N(w)\backslash\{u,v\}$. If
$N(w)\backslash\{u,v\}$ contains an induced $P_3$, then $\langle\{w,V(P_3)\}\rangle \cong diamond$, a contradiction. Hence $G$ is $(P_3 \cup P_2)$-free.  Also, $\omega(G)=3$ and $\chi(G)=6$ (see, \textcolor{blue}{www.win.tue.nl/\textasciitilde aeb/graphs/Schlaefli.html}). In addition, the Gr\"{o}tzsch graph $\mu(C_5)$ given in Figure \ref{gro} is an example of a $\{P_3\cup P_2, diamond\}$-free graph with $\omega(\mu(C_5))=2$ and $\chi(\mu(C_5))=4$. Thus the bounds given in Theorem \ref{diamond} are optimal for all the values of $\omega(G)$.

\subsection*{Acknowledgment}
\small The first author's research was supported by the Council of Scientific and Industrial Research,  Government of India, File No: 09/559(0133)/2019-EMR-I.
The second author's research was supported by Post Doctoral Fellowship at Indian Institute of Technology, Palakkad.

\renewcommand{\baselinestretch}{1}\small
\bibliographystyle{ams}
\bibliography{chi_Binding_Functions_for_classes_P_3cupP_2_free_Graphs}

\end{titlepage}
\end{document}